\documentclass[11pt,twoside]{article}
\usepackage{bbm}
\usepackage{amsmath, amssymb, mathrsfs, graphicx}

\def\scr{\mathscr}

\def\gz{\gamma}  
\def\lz{\lambda} 
     \def\oz{\omega}

        \def\uz{\theta}

 \def\ddz{\Delta}
    
  \def\ooz{\Omega}

\def\qd{\quad}
\def\qqd{\qquad}

\newcommand{\mathsym}[1]{{}}

\def\scr{\mathscr}

\def\le{\leqslant}
\def\ge{\geqslant}

\font\cms=cmss9 scaled \magstep1

\def\nnd{\noindent}

\def\thm{\nnd\bg{thm1}}
\def\crl{\nnd\bg{crl1}}
\def\lmm{\nnd\bg{lmm1}}
\def\prp{\nnd\bg{prp1}}

\def\xmp{\nnd\bg{xmp1}}

\def\rmk{\nnd\bg{rmk1}}

\def\dethm{\end{thm1}}
\def\decrl{\end{crl1}}
\def\delmm{\end{lmm1}}
\def\deprp{\end{prp1}}
\def\dexmp{\end{xmp1}}

\def\dermk{\end{rmk1}}

\def\prf{\medskip \noindent {\bf Proof}. }
\def\qed{\text{\quad $\square$}}
\def\deprf{\qed\medskip}
\def\bg{\begin}
\def\be{\bg{equation}}
\def\de{\end{equation}}

\def\dear{\end{eqnarray}}
\def\lb{\label}
\def\ct{\cite}

\newcommand{\rf}[2]{[\ref{#1}; #2]}

\def\den{\end{enumerate}}

\def\d{\text{\rm d}}

\allowdisplaybreaks[4]

\begin{document}
%\song

%\setcounter{page}{1016}
\allowdisplaybreaks[4]

\renewcommand{\thefootnote}{\fnsymbol{footnote}}

\noindent {Commun. Math. Stat. (2014) 2:17--32}

\thispagestyle{empty}
\renewcommand{\thefootnote}{\fnsymbol{footnote}}

%\noindent {Acta Mathematica Sinica, English Series}\newline
%\noindent {January 2013, Volume 29, Issue 1, pp 1-32 }

\vspace*{.5in}
\begin{center}
{\bf\Large Isospectral operators}
\vskip.15in {Mu-Fa Chen\qqd \qqd\qqd\qqd \qqd Xu Zhang}
\end{center}
\begin{center} (Beijing Normal University) \qqd (Beijing University of Technology)\\
\vskip.1in March 6, 2014
%\vskip.1in November 17, 2013
\end{center}
\vskip.1in

\markboth{\sc Mu-Fa Chen and Xu Zhang}{\sc Isospectral operators}

%\title{Exponential Convergence Rate in Entropy}

%\author{Mu-Fa Chen}

\date{}
%{February 24, 2013}

%\maketitle

%\footnotetext{Received May 17, 2012; accepted June 18, 1012}
%\footnotetext{2000 {\it Mathematics Subject Classifications}.\quad 26D10, 60J60, 34L15.}
%\footnotetext{{\it Key words and phases}.\quad
%Hardy-type inequality, optimal constant, variational formulas, approximating procedure.}

\begin{abstract}

For a large class of integral operators or second order differential operators,  their isospectral (or cospectral)
operators are constructed explicitly in terms of $h$-transform (duality). This provides us a simple way to
extend the known knowledge on the spectrum
(or the estimation of the principal eigenvalue) from a smaller class of operators to a much larger one.
In particular, an open problem about the positivity of the principal eigenvalue for birth--death processes is solved in the paper.
\end{abstract}

\nnd {\small 2000 {\it Mathematics Subject Classification}: 58J53; 37A30.}

\nnd {\small {\it Key words and phases}. Isospectral; harmonic function; integral operator, differential operator.}

\section{Introduction}

Let us consider the elliptic operators
$$\aligned
L&= \sum_{i, j}a_{ij}(x)\partial_{ij}^2+ \sum_i b_i(x)\partial_i +  c(x),\\
\widetilde L&= \sum_{i, j}{\tilde a}_{ij}(x)\partial_{ij}^2+ \sum_i {\tilde b}_i(x)\partial_i
\endaligned$$
on $L^2(\mu)$ and $L^2(\tilde\mu)$ (real) respectively, where
$\tilde\mu=h^2\mu$ for a given measure $\mu$ and some $h\ne 0$.
Their main difference is that $c(x)\not\equiv 0$. We are
interested in when the operators $L$  and $\widetilde L$ are
$L^2$-isospectral in the following sense
$$(L f, f)_{\mu}=\big({\widetilde L}{\tilde f}, {\tilde f}\big)_{\tilde\mu},
\qqd \text{for every }{\tilde f}:=f/h, \; f\in {\scr D}(L).$$
Here is one of our typical results in the note (cf. Theorems \ref{t3-1} and \ref{t3-6} in Section 3).

\thm\lb{t1-1}
{\cms
\begin{itemize}
\item[(1)] Given $L$ on $L^2(\mu)$ having domain ${\scr D}(L)$, let $h\ne 0$, $\mu$-a.e. be $L$-harmonic: $Lh=0$, $\mu$-a.e.,
   then $L$ is $L^2$-isospectral to $\widetilde L$:
   $$\aligned
   \widetilde L& =L_0 + 2 h^{-1}\langle a \nabla h, \nabla \rangle,\qqd
   {\scr D}\big(\widetilde L\big)=\{f: fh\in {\scr D}(L)\}.
   \endaligned$$
   where $L_0=L-c$.
\item[(2)] Given $\widetilde L$ on $L^2(\tilde\mu)$ having domain ${\scr D}\big(\widetilde L\big)$,
  then for each $h\ne 0$, $\mu$-a.e., $\widetilde L$ is $L^2$-isospectral to $L$:
$$\aligned
&L=\widetilde L- \frac{2}{h}\big\langle {\tilde a} \nabla h, \nabla \big\rangle
  +\bigg[\frac{2}{h^2}\big\langle {\tilde a} \nabla h, \nabla h \big\rangle
       -\frac{1}{h}{\widetilde L}h \bigg],\\
      & {\scr D}(L)=\big\{f: f/h\in {\scr D}\big(\widetilde L\big)\big\},
       \endaligned$$
       where $\langle \cdot, \cdot \rangle$ denotes the Euclidean inner product.
\end{itemize}
}
\dethm

As a typical application of Theorem \ref{t1-1}, we obtain the next result.
To state it, we need to explain the meaning of eigenvalue in different sense. We say that $\lz$ is an
eigenvalue of $L$ in the ordinary sense if $L g= \lz g$ for some $g\ne 0$. It is called
a $L^2$-eigenvalue if additionally, $g\in L^2(\mu)$.

\crl\lb{t1-2}{\cms For each $h\in {\scr C}^2(\mathbb R)$, $h\ne 0$, a.e., the operator
$$\aligned
&L^h =\frac 1 2 \frac{\d^2}{\d x^2 } -\bigg(x+ \frac{h'}{h}\bigg) \frac{\d }{\d x}
 +\bigg[\bigg(\frac{h'}{h}\bigg)^2+ x \frac{h'}{h} -\frac{h''}{2h}\bigg]
\endaligned$$
has $L^2$-eigenvalues $\lz_n\big(L^h\big)=-n$ with eigenfunctions
$$g_n(x)=(-1)^n h(x)e^{x^2}\frac{\d^n}{\d x^n}\big(e^{-x^2}\big),\qqd n\ge 0,$$
respectively.  A particular class of $L^h$ is the following
$$L^b=\frac 1 2 \frac{\d^2}{\d x^2 } -b(x) \frac{\d }{\d x}
+\frac 1 2 \big[b(x)^2-b'(x)-x^2+1\big],\qqd b\in {\scr C}^1(\mathbb R).$$}
\decrl

\prf Noting that the Ornstein-Uhlenbeck operator
$$\aligned
& {\widetilde L}=\frac 1 2 \frac{\d^2}{\d x^2 } -x \frac{\d }{\d x},\qqd
 {\scr D}\big({\widetilde L}\big)\supset {\scr C}_0^{\infty}(\mathbb R)
\endaligned
$$
has ordinary eigenvalues $\lz_n\big(\widetilde L\big)=-n$ with eigenfunctions
$$g_n(x)=(-1)^n e^{x^2}\frac{\d^n}{\d x^n}\big(e^{-x^2}\big),\qqd n\ge 0,$$
respectively (cf. \rf{cmf12a}{Example 5.1}). Clearly,
the polynomial function $g_n\in L^2(\tilde\mu)$ for every $n\ge 0$, where
$\tilde\mu(\d x)=\exp(-x^2)\d x$.
Hence, the eigenvalues are all $L^2$-ones.
Now, the first assertion follows from part (2) of Theorem \ref{t1-1}.
The last assertion then follows by setting $h=\exp\psi$ with $\psi'=b-x$:
$$\bigg(\frac{h'}{h}\bigg)^2\!+ x \frac{h'}{h} -\frac{h''}{2h}
={\psi'}^2\!+x\psi' -\frac 1 2 \big(\psi''+{\psi'}^2\big)
= \psi'\bigg(x+\frac 1 2 \psi'\bigg)-\frac 1 2 \psi''.\qqd \square$$

Corollary \ref{t1-2} says that a large class of operators are all
isospectral to the rather simple Ornstein-Uhlenbeck operator. This
indicates the value of the study on isospectral operators. It
should be pointed out that the technique is still valuable even if
you know only some estimates of the principal eigenvalue of
$\widetilde L$ but have no knowledge on the other part of the
spectrum of $\widetilde L$, since our knowledge on the principal
eigenvalue of $L$ is still rather limited.

Actually, Theorem \ref{t1-1} comes from a very simple observation. For completeness, here we
write its complex version, even though we will use only its real version later on.

\lmm\lb{t1-3}{\cms Let $(E, {\scr E}, \mu)$ be a measure space and let $h$ be Lebesgue
measurable: $E\to {\mathbb C}$, $h\ne 0$, $\mu$-a.s. Then
\bg{itemize}
\item[(1)] ${\tilde f}:=\mathbbm{1}_{[h\ne 0]}f/h$ is an isometry from $L^2(E, \mu)$ to
  $L^2(E, \tilde\mu)$ (complex), where $\tilde\mu=|h|^2 \mu$.
\item[(2)] Let $L$ be an operator on $L^2(E, \mu)$ with domain ${\scr D}(L)$.
Define an operator $\widetilde L$ as follows:
\be {\widetilde L}{\tilde f}= \mathbbm{1}_{[h\ne 0]}\frac{1}{h} L \big({\tilde f}h\big),
\qqd {\scr D}\big(\widetilde L\big)=\big\{{\tilde f}\in {\scr E}:
{\tilde f}h \in {\scr D}(L)\big\}.
\lb{1-3}\de
Then the operators $(L, {\scr D}(L))$ on $L^2(E, \mu)$
and $\big(\widetilde L, {\scr D}\big(\widetilde L\big)\big)$ on $L^2\big(E, \tilde\mu\big)$ are isospectral
(say $L$ and $\widetilde L$ are $L^2$-isospectral, for short) (in the following sense):
$$(Lf, f)_{\mu}= \big({\widetilde L}{\tilde f}, {\tilde f}\big)_{\tilde\mu},\qqd f\in {\scr D}(L).$$
\item[(3)] If additionally, $h\in {\scr D}(L)$, then $\widetilde L \mathbbm{1}=0$, ${\tilde\mu}$-a.e. iff $h$ is
$L$-harmonic: $Lh=0$, $\mu$-a.s.
\end{itemize}}
\delmm

\prf Recall the inner product in a complex $L^2$-space:
$$(f, g)_{\mu}=\int_E f {\bar g}\d\mu.$$
The first assertion is obvious:
$$\int_E |f|^2\d\mu =\int_{E\, [h\ne 0]}  \big|\tilde f\big|^2 |h|^2\d\mu= \int_E \big|\tilde f\big|^2\d{\tilde\mu}.$$
By definition, for $\tilde f\in {\scr D}\big(\widetilde L\big)$, we have
$\tilde f h\in {\scr D}(L)\subset L^2(E, \mu)$. Then we have not only
${\tilde f}\in L^2\big(E, \tilde \mu\big)$ but also $L\big(\tilde f h\big)\in L^2(E, \mu)$.
This means that ${\widetilde L}{\tilde f} \in L^2\big(E, \tilde\mu\big)$. Hence, as an operator on
$L^2\big(E, \tilde\mu\big)$, ${\widetilde L}$ is well defined. Furthermore, we have
$$(Lf, f)_{\mu}\!=\!\big({L}\big({\tilde f}h\big), {\tilde f}{h}\big)_{\mu}
\!=\!\int_E {\overline{{\tilde f}h}}\, {L}\big({\tilde f}h\big) \d\mu
\!=\!\int_E {\bar{\tilde f}} ({\bar h}h) \frac{1}{h}{L}\big({\tilde f}h\big) \d\mu
\!=\!\big({\widetilde L}{\tilde f}, {\tilde f}\big)_{\tilde\mu}. $$
We have thus proved the second assertion.
Clearly, if $h\in {\scr D}(L)$, then $\mathbbm{1} h=h\in L^2(E, \mu)$ and hence
$\mathbbm{1} \in L^2\big(E, \tilde\mu\big)$ which implies that $\tilde\mu (E)<\infty$. Furthermore,
$\mathbbm{1}\in {\scr D}\big(\widetilde L\big)$ by definition of ${\scr D}\big(\widetilde L\big)$.
Therefore, the last assertion follows by definition of $\widetilde L$.
\deprf

For non-symmetric operators, their spectrum can be complex. Hence, it is natural to use the
complex $L^2$-theory. However, in this note, we use the real $L^2$-spaces only. Thus,
the $L^2$-isospectral (real) here means the spectrum of their symmetrized operators.
The last assertion of the lemma suggests us, as we will do often later, to choose $h$ as an
$L$-harmonic function in a weak (pointwise) sense (in other words, $h$ is in a weak domain of $L$)
without assuming $h\in {\scr D}(L)$.
Then $\widetilde L \mathbbm{1}=0$ is meaningful in the weak sense. In this way, we can
construct the operator $\widetilde L$ explicitly, which is the main goal of this note.
Furthermore, part (3) of the lemma has the following extension.
\rmk\lb{t1-4}{\cms For fixed $B\in {\scr E}$, $\widetilde L \mathbbm{1}=0$, ${\tilde\mu}$-a.e. on $B$ iff
$Lh=0$, $\mu$-a.s. on $B$.}\dermk

We will illustrate later an application of this assertion in the context of Markov chains.
Clearly, the $L$-harmonic
function is an eigenfunction corresponding to the eigenvalue $\lz=0$. However, $\lz=0$ is not
necessary an eigenvalue in the $L^2$-sense unless $h\in L^2(E, \mu)$.

One may write $\widetilde L= h^{-1} L(h\, \bullet)$ ($\mu$-a.e.) for short. Because of this,
$\widetilde L$
is called a $h$-transform of $L$. Alternatively, define an operator $H$:
$$Hf= hf, \qqd {\scr D}(H)=\{f\in L^2(E,\mu): hf\in {\scr D}(L)\}.$$
Then, we indeed have $\widetilde L=H^{-1} L H$. In view of this, $L$ and $\widetilde L$
are similar and so are $L^2$-isospectral. More generally (without assuming the invertibility of $H$),
$$H \widetilde L = L H.$$
Because of this, $L$ and $\widetilde L$ are called dual with respect to $H$.
Therefore, the $h$-transform is indeed a special duality. For a different dual,
refer to \rf{cmf10}{\S 5 and \S 10}. Note that in the later case, we were
interested in the principal eigenvalue only, but the transform used there is still isospectral.
The reason is that the isospectral transform is easier to handle even though it
looks rather strong. We remark that when
$E$ has boundary $\partial E$, one may deduce a boundary condition for $\widetilde L$ from
that of $L$, based on the transform $\tilde f=\mathbbm{1}_{[h\ne 0]} f/h$.

Having figured out the dual operators, in the study of their spectrum for Markov processes, it is more convenient in practice to use their
extension to the Dirichlet forms, especially for the operator $\big(\widetilde L, {\scr D}\big(\widetilde L\big)\big)$. Generally speaking, Lemma \ref{t1-3} says that for a given
Dirichlet form $(D, {\scr D}(D))$ on $L^2(\mu)$, its dual form $\big({\widetilde D}, {\scr D}\big({\widetilde D}\big)\big)$ on $L^2(\tilde\mu)$ is given by
$${\widetilde D}\big({\tilde f}\big)=D\big(\tilde f h, \tilde f h\big),\qqd
{\scr D}\big({\widetilde D}\big)=\big\{\tilde f\in {\scr E}: {\tilde f}h \in {\scr D}(D)\big\}.$$
Certainly, one may go to the inverse way, defining $(D, {\scr D}(D))$ in terms of
$\big({\widetilde D}, {\scr D}\big({\widetilde D}\big)\big)$.
In particular, for the O.-U. operator
used in the proof of Corollary \ref{t1-2}, corresponding to
$\big(\widetilde L, {\scr D}\big(\widetilde L\big)\big)$,
the Dirichlet form $\big({\widetilde D}(f), {\scr D}\big(\widetilde D\big)\big)$ is
$$\aligned
&{\widetilde D}(f)= \int_{\mathbb R} {f'}^2e^{-x^2}\d x,\\
&{\scr D}\big(\widetilde D\big)\!=\!\big\{f\in L^2(\tilde\mu): {\widetilde D}(f)\!<\!\infty\big\}
\!=\!\bigg\{f: \int_{\mathbb R} \big[f^2+  {f^{\prime}}^2 \big]e^{-x^2}\d x\!<\!\infty\bigg\}.
\endaligned$$
In the case that the potential term $c^h$ (the last term) in $L^h$ is non-positive, then
$L^h$ corresponds to the operator of a diffusion having killing rate $-c^h$, to which we certainly
have a Dirichlet form $\big(D^h, {\scr D}\big(D^h\big)\big)$ on $L^2(\mu^h)$:
$$\aligned
&D^h(f)=\int_{\mathbb R}\big[ {f'}^2(x)-c^h (x){f}^2(x)\big]e^{-x^2}\frac{\d x}{h(x)^2},\\
&{\scr D}\big(D^h\big)\!=\!\bigg\{f:
\int_{\mathbb R} \big[f^2+  (f' h -f h')^2 \big]e^{-x^2}\d x\!<\!\infty\bigg\},\\
&c^h(x)= \bigg[\bigg(\frac{h'}{h}\bigg)^2+ x \frac{h'}{h} -\frac{h''}{2h}\bigg](x),
\qqd \mu^h(\d x)= e^{-x^2} \frac{\d x}{h(x)^2}.
\endaligned$$
Here ${\scr D}\big(D^h\big)$ is deduced from ${\scr D}\big(\widetilde D\big)$,
based on Lemma \ref{t1-3}.
For general $c^h(x)\in {\mathbb R}$, this symmetric form may not be a Dirichlet one
even though it does have nonnegative spectrum in view of our isospectral property.
Actually, Lemma \ref{t1-3} is meaningful in a very general setup rather than Markov
processes.

The $h$-transform, or the Doob's $h$-transform is a well-known topic in probability/potential
theory. Here we mention only two related papers \ct{prg09, wj12} where the tool is used to
study the principal eigenvalue. In \ct{prg09}, the following model
$$\aligned
&L=\frac 1 2 \frac{\d}{\d x} a \frac{\d}{\d x} - \frac 1 2 \bigg(\frac{b^2}{a}+b'\bigg),\\
&{\widetilde L}=\frac 1 2 \frac{\d}{\d x} a \frac{\d}{\d x}+ b \frac{\d}{\d x},\\
&h(x)=\exp\bigg[\int_0^x \frac{b}{a}(y)\d y\bigg]
\endaligned
$$
is carefully handled and applied to multi-dimensional diffusion operators.
In \ct{wj12}, a class of symmetric Markov processes having killings are studied
and some upper and lower estimates for the first eigenvalue are presented.

The remainder of this note is organized as follows. In the next two sections,
we apply Lemma \ref{t1-3}, respectively, to two special classes of operators: either integral operators
for Markov pure jump processes or the operators for diffusions.

\section{Integral operators}

\thm\lb{t2-1}{\cms Let $(q(x), q(x, \d y))$ be a totally stable and conservative $q$-pair on $(E, {\scr E}, \mu)$
(cf. \rf{cmf04}{Definition 1.9}).
For a given function $c\in {\scr E}$ with $c\le q$, define an operator $\ooz$
$$\ooz f(x)=\int_E q(x, \d y)\big[f(y)-f(x)\big]+c(x)f(x),\qqd x\in E$$
with domain ${\scr D}(\ooz)\subset L^2(E, \mu)$.
Next, let $h\,(> 0, \mu\text{\cms -a.e.})$ be $\ooz$-harmonic (if exists): $\ooz\, h=0$, $\mu$-a.e.
on $E.$
Define a new totally stable and conservative $q$-pair $\big(\tilde q(x), \tilde q(x, \d y)\big)$ as follows.
$$\aligned
{\tilde q}(x, A)&=  \mathbbm{1}_{[h(x)\ne 0]}\frac{1}{h(x)}\int_A q(x, \d y) h(y),\qqd A\in {\scr E},\\
 {\tilde q}(x)&={\tilde q}(x, E),
 \qd \mu\text{\cms -a.e. }x\in E. \endaligned$$
 Set
$$\aligned
 {\widetilde \ooz}f(x)&=\int_E{\tilde q}(x, \d y)\big[f(y)-f(x)\big],\qqd \mu\text{\cms -a.e. }x \in E,\\
 {\scr D}\big(\widetilde \ooz\big)&=\big\{{\tilde f}\in {\scr E}:
{\tilde f}h \in {\scr D}(\ooz)\big\}.
\endaligned$$
Then $\ooz$ and $\widetilde \ooz $ are $L^2$-isospectral.
}
\dethm

\prf Noting that $h\,(>0, \mu$-a.e.) is $\ooz$-harmonic by assumption, we have
$$\big[q(x)-c(x)\big]h(x)=\int_E q(x, \d y) h(y)\ge 0.$$
Hence $h$ is $q(x, \cdot)$-integrable for a.e.-$x\in E$ and moreover $q\ge c$.
Therefore, the new $q$-pair $\big(\tilde q(x), \tilde q(x, \d y)\big)$ is totally
stable. It is clearly conservative.
By definition of $\widetilde \ooz$, we have on the set $[h> 0]$,
$$\aligned
{\widetilde \ooz} (f)(x)&=\int_E{\tilde q}(x, \d y)\big[f(y)-f(x)\big]\\
&=\frac{1}{h(x)} \int_E{q}(x, \d y)\big\{\big[(fh)(y)-(fh)(x)\big]+f(x)\big[h(x)-h(y)\big]\big\}\bigg]\\
&=\frac{1}{h(x)}\bigg[\int_E\!{q}(x, \d y)\big[(fh)(y)\!-\!(fh)(x)\big]
\!-\!\!f(x)\!\! \int_E\!{q}(x, \d y)\big[h(y)\!-\!h(x)\big]\\
&=\frac{1}{h(x)}\big[\ooz (fh)(x)-c(fh)(x)-f(x)[\ooz h(x)-(ch)(x)]\big]\\
&=\frac{1}{h(x)}\big[\ooz (fh)(x)-f(x)\ooz h(x)\big].
\endaligned
$$
Now, by harmonic property of $h$, the right-hand side is equal to
$$\frac{1}{h(x)}\ooz (fh)(x)\qqd\text{ on } [h> 0].$$
The assertion then follows from Lemma \ref{t1-3}.
\deprf

We mention that the positive condition of $h$ used in the theorem is to
keep $\big(\tilde q(x), \tilde q(x, \d y)\big)$ to be a $q$-pair. This
is certainly not necessary in a general context: considering general
integral kernel instead of the nonnegative one.

The inverse of the last theorem goes as follows.

\thm{\cms Given a totally stable and conservative $q$-pair $\big(\tilde q(x), \tilde q(x, \d y)\big)$
and a positive ${\scr E}$-measurable function $h$ such that $h^{-1}$ is
${\tilde q}(x, \cdot)$-integrable for each $x\in E$, the operator
$\big(\widetilde \ooz, {\scr D}\big(\widetilde \ooz\big)\big)$ on $L^2(E, \tilde\mu)$ corresponding to
the $q$-pair $\big(\tilde q(x), \tilde q(x, \d y)\big)$ is $L^2$-isospectral to the following operator
$\ooz$ on $L^2(E, \mu)$ \big($\mu:= h^{-2}{\tilde\mu}$\big):
$$\aligned
\ooz f(x)&=\int_E q(x, \d y)[f(y)-f(x)]+ c(x)f(x),\\
{\scr D}(\ooz)&=\big\{f\in {\scr E}: f/h\in {\scr D}\big(\widetilde \ooz\big)\big\}\subset L^2(E, \mu),
\endaligned$$
where
\begin{gather}
q(x, \d y)= h(x) \frac{{\tilde q}(x, \d y)}{h(y)},\nonumber\\
c(x)= \int_E {\tilde q}(x, \d y)\bigg[\frac{h(x)}{h(y)}-1\bigg],\qqd x\in E.\nonumber
\end{gather}
}
\dethm

\prf It is simply a use of the duality $\ooz=H {\widetilde \ooz}H^{-1}$,
noting the property that $\ooz h=0$ is now automatic since $\widetilde \ooz 1=0$.
The remainder of the proof is mainly a careful computation.
\deprf

It is the place to discuss the existence of a positive $\ooz$-harmonic function.
Let $c(x)<q(x),\; x\in E$. Choose and fix a reference point $\uz\in E$.  By
\rf{cmf04}{Theorem 2.2}, there exists uniquely the minimal solution $(h^*(x): x\in E)$
with $h^*(\uz)=1$ to the following nonnegative equation
\be h(x)= \int_{E\setminus \{\uz\}}\frac{q(x, \d y)}{q(x)-c(x)}h(y)+
\frac{q(x, \{\uz\})}{q(x)-c(x)},\qqd x\ne \uz.\lb{2-5}\de
Moreover, the solution can be obtained in the following way: let
$$\aligned
h^{(1)}(x)&= \frac{q(x, \{\uz\})}{q(x)-c(x)},\qqd x\ne \uz,\\
h^{(n+1)}(x)&= \int_{E\setminus \{\uz\}}\frac{q(x, \d y)}{q(x)-c(x)}h^{(n)}(y)+
\frac{q(x, \{\uz\})}{q(x)-c(x)},\qqd x\ne \uz,\; n\ge 1.
\endaligned$$
Then  for each $x\ne \uz$, $h^{(n)}(x)\uparrow h^*(x)\in [0, \infty]$ as $n\to\infty$.

\prp\lb{t2-3}{\cms Let $c(x)<q(x)$ for every $x\in E$ and assume that $q(x, \{\uz\})>0$ for some $x\ne\uz$. Then
the equation $\ooz h=0$ has a non-trivial (finite) solution iff the minimal solution
$(h^*(x): x\in E)$ to (\ref{2-5}) is finite. Equivalently, there is a finite $f$
satisfying the inequality
$$f(x)\ge \int_{E\setminus \{\uz\}}\frac{q(x, \d y)}{q(x)-c(x)}f(y)+
\frac{q(x, \{\uz\})}{q(x)-c(x)},\qqd x\ne \uz.$$
Then we actually have $f(x)\ge h^*(x)$ for every $x\in E$.\footnote{\rm Correction. Here the uniqueness of the solution $h$
to the equation $\Omega h=0$ with
$h(\theta)=1$ up to a positive constant is needed. Otherwise,
$(h^*(x): x\in E)$ is only a lower
bound of $h$.}}
\deprp

\prf For a given finite non-trivial $\ooz$-harmonic function $h$,
choosing $h(\uz)=1$,
one may write down immediately equation (\ref{2-5}).

Conversely, a finite solution $h^*$ to (\ref{2-5}) is clearly a
$\ooz$-harmonic function.
From the construction given above, it is also clear that $h^*(x)>0$ once
$q(x, \{\uz\})>0$. The last assertion of the proposition is essentially a comparison theorem
\rf{cmf04}{Theorem 2.6}.
\deprf

It is clear from the proof above, to obtain a positive harmonic $h$, some irreducible
condition is necessary. Noting that it is often practical to find an explicit comparison function
$f$, and $h^{(n)}$ for each $n$ is already explicit, we have explicit estimates of $h^*$
which may not be easy to obtain explicitly.

Before moving further, we discuss an alternative way to describe the $\ooz$-harmonic
function. Suppose that $\sup_x c(x)<\infty$. Then by a shift if
necessary, we may and will assume for a moment that $\sup_x c(x)\le 0$. Define
$$\aligned
z^{(0)}(x)&= 1,\qqd x\in E,\\
z^{(n+1)}(x)&= \int_{E}\frac{q(x, \d y)}{q(x)-c(x)}z^{(n)}(y),\qqd x\in E,\; n\ge 1.
\endaligned$$
Then $z^{(n)}(x)\downarrow {\bar z}(x)$ as $n\to\infty$ for each $x\in E$.
This is an analog of the maximal exit solution in the study of $q$-processes,
cf. \rf{cmf04}{Lemma 2.39}. The proof for the conclusion is easy, simply use
the property
$$\frac{q(x, E)}{q(x)-c(x)}\le 1,\qqd x\in E.$$

\rmk{\cms Let $\sup_x c(x)\le 0$. Then a bounded $\ooz$-harmonic function is non-zero iff so is the maximal solution
${\bar z}$ constructed above.
}
\dermk

To apply the previous results, Theorem \ref{t2-1} for instance, to finite state spaces, say $E=\{0, 1, \ldots, N\}$ for some
$N\ge 3$, one meets a problem about the existence of positive $\ooz$-harmonic $h$. For which,
there $N+1$ homogeneous equations with $N+1$ variables $h_0, h_1, \ldots, h_N$. Because of
the homogeneous property in $h$, one may assume that $h_0=1$ once a non-trivial solution $h$ exists
with $h_0\ne 0$ for instance. Thus, we have only $N$ free variables in $N+1$ equations. Then a finite
non-trivial solution often does not exist (or equivalently, the minimal solution given in Proposition \ref{t2-3}
may be infinite). To overcome this difficulty, one has to decrease the number of
equations. This is the reason we will adopt a local harmonic condition below. Then, one
needs non-trivial ${\tilde c}_i$ in the corresponding operator $\widetilde \ooz$.

\thm{\cms Let $E=\{0, 1, \ldots, N\}$ for some $N\ge 3$ and let $Q=(q_{ij})$ be a conservative
$Q$-matrix on $E$. For given $(c_i: i=0, 1, \ldots, N)$ with $c_i\le q_i:=-q_{ii}$ for $i=0, 1, \ldots, N-1$,
set $\ooz=Q+\text{\rm diag}(c_i)$. Next, let $h> 0$ be $\ooz$-harmonic on
$\{0, 1, \ldots, N-1\}$, i.e.,
$$\ooz\, h=0 \qqd \text{\cms on } \{0, 1, \ldots, N-1\}.$$
Define ${\tilde q}_{ij}\; (i, j\in E)$ as in Theorem \ref{t2-1}:
$${\tilde q}_{ij}=h_i^{-1} q_{ij}h_j, \qqd i, j\in E.$$
 Next, define
$\tilde c_i=0$ on $\{0, 1, \ldots, N-1\}$ and
$$\tilde c_N^{}=c_N^{}+ \sum_{j\le N} q_{Nj}^{} \bigg(\frac{h_j}{h_N}-1\bigg).$$
Denote by $\widetilde \ooz$ the operator corresponding to the matrix
$\big(\tilde q_{ij}\big)+\text{\rm diag}(\tilde c_i)$.
Then $\ooz$ and $\widetilde \ooz$ are $L^2$-isospectral.}
\dethm

\prf Following the proof of Theorem \ref{t2-1}, restricted to $\{0, 1, \ldots, N-1\}$, we see that
$${\widetilde \ooz}{\tilde f}(i)=\frac{1}{h_i} \ooz\big(\tilde f h\big)(i) \qqd \text{on }\{0, 1, \ldots, N-1\}.$$
We now show that this equality also holds for $i=N$.
$$\aligned
{\widetilde \ooz} f(N)&= \sum_{j\le N} {\tilde q}_{Nj}^{} (f_j-f_N)+{\tilde c}_N f_N\\
&=\frac{1}{h_N} \sum_{j\le N} q_{Nj}^{} \big[(fh)_j-(fh)_N\big]-\frac{f_N}{h_N}
\sum_{j\le N} q_{Nj}^{} (h_j-h_N) +{\tilde c}_N f_N\\
&=\frac{1}{h_N} Q (fh)(N) -\frac{1}{h_N} c_N h_N f_N-\frac{f_N}{h_N}
\sum_{j\le N} q_{Nj}^{} (h_j-h_N) +{\tilde c}_N f_N\\
&=\frac{1}{h_N} \ooz (fh)(N).
\endaligned
$$
From Remark \ref{t1-4}, it follows that $c_i=0$ on $\{0, 1, \ldots, N-1\}$.
The required assertion now follows from Lemma \ref{t1-3}.\deprf

A typical application of Theorem \ref{t2-1} to the single birth
processes is presented in \ct{zx13}. In this case, the
$\ooz$-harmonic function has a very simple expression (cf. \rf{chzh}{Theorem 1.1}).
In particular, for the killing case, the function is not only positive but also
non-decreasing. It is interesting to note that for single birth
processes, the function $h$-dual is again the same type, but the
measure $\mu$-dual
$${\bar q}_{ij}=\frac{\mu_j q_{ji}}{\mu_i}, \qqd i, j\in E$$
maps the single birth type to the single death type.
Next, for birth--death processes
with birth and death rates $b_i$ and $a_i$, respectively, and with killing rates $-c_i\ge 0$,
we have
$${\tilde a}_i = a_i \frac{h_{i-1}}{h_i}\,(\le a_i),\;i\ge 1, h_0=1,\qqd
{\tilde b}_i = b_i \frac{h_{i+1}}{h_i}\,(\ge b_i), \; i\ge 0.$$
Then
$${\tilde\mu}_i=\frac{{\tilde b}_0\ldots {\tilde b}_{i-1}}{{\tilde a}_1\ldots {\tilde a}_{i}}
=\frac{{b}_0\ldots {b}_{i-1}}{{a}_1\ldots {a}_{i}}h_i^2=h_i^2\mu_i,\qqd
{\hat{\tilde\nu}}_i=\frac{1}{{\tilde\mu}_i{\tilde b}_i}= \frac{1}{h_i h_{i+1}}{\hat\nu}_i,\qqd i\ge 0.$$
For finite state space, we have
$${\tilde c}_N^{}=c_N^{}+ a_N^{}\bigg(\frac{h_{N-1}}{h_N}-1\bigg).$$
Clearly, ${\tilde c}_N^{}\le 0$ since so does $c_N^{}$. However, the story is still meaningful
for general $c_i\in {\mathbb R}$ satisfying $c_i\le a_i+b_i$ for all $i\ge 0$.

To conclude this section, we answer an open question for birth--death processes with state
space $\{0, 1, 2, \ldots\}$. For this,
we need some notation. Given birth rates $b_i>0$\,($i\ge 0$), death rates $a_i>0$\,($i\ge 1$) and killing rates $-c_i\ge 0$\,($i\ge 0$), define
\begin{gather}
{\tilde q}_n^{(k)}=
{\begin{cases}
-c_n, \quad & 0\le k\le n-2\\
a_n-c_n,  & k=n-1,
\end{cases}}\nonumber\\
{\widetilde F}_i^{(i)}=1,\;\;
{\widetilde F}_n^{(i)}=\frac{1}{b_n}\sum_{k=i}^{n-1}{\tilde q}_n^{(k)} {\widetilde F}_k^{(i)},
\qquad n>i\ge 0,\nonumber\\
h_n=1-\sum_{0\le k\le n-1} \sum_{0\le j\le k} {\widetilde F}_k^{(j)}\frac{c_j}{b_j},\qquad n\ge 0.\nonumber
\end{gather}
Next, define the principal eigenvalue $\lambda_0$ as follows.
$$\lambda_0=\inf\bigg\{\sum_{k\ge 0}\mu_k\big[b_k (f_{k+1}-f_k)^2-c_k f_k^2\big]:
\sum_{k\ge 0} \mu_k f_k^2=1,\; f\text{ has finite support}\bigg\}.$$
Here is a solution to the Open Problem 9.13 in \cite{cmf10}.

\thm{\cms For birth--death processes as above, we have ${\tilde \delta}\le \lambda_0^{-1}\le 4 {\tilde \delta}$, where
$${\tilde \delta}=\sup_{n\ge 0}\sum_{j=0}^n {\tilde \mu}_j\sum_{k\ge n}\hat{\tilde \nu}_k
=\sup_{n\ge 0}\sum_{j=0}^n {\mu}_jh_j^2\sum_{k\ge n} \frac{1}{h_k h_{k+1}\mu_k b_k}.$$
In particular, $\lambda_0>0$ iff ${\tilde \delta}<\infty$.}
\dethm

\prf The harmonic function $h$ we need for applying Theorem \ref{t2-1} is given by
\rf{chzh}{Theorem 1.1}. Then the result follows by applying \rf{cmf10}{Theorem 3.1}
to the process with rates $(\tilde b_i, \tilde a_i)$ and using $\tilde \mu_i$ and
$\hat{\tilde \nu}_k$ just computed above.\deprf
\section{Differential operators}

We now turn to study the second-order differential operators.

\thm\lb{t3-1}{\cms Consider the elliptic operator
$$L= \sum_{i, j}a_{ij}(x)\partial_{ij}^2+ \sum_i b_i(x)\partial_i +  c(x)$$
with a domain ${\scr D}(L)$, and let $h\ne 0$ a.e. (with respect to Lebesgue measure) be $L$-harmonic. Here
$$\partial_i=\d/\d x_i,\qqd \partial_{ij}^2=\partial_i\partial_j.$$
Define
$$\widetilde L= \sum_{i, j}{\tilde a}_{ij}(x)\partial_{ij}^2+ \sum_i {\tilde b}_i(x)\partial_i, $$
with domain ${\scr D}\big(\widetilde L\big)$ defined in Lemma \ref{t1-3}, where
$${\tilde a}_{ij}(x)={a}_{ij}(x),\qqd
{\tilde b}_i(x)= b_i(x)+ \frac{2}{h(x)} \sum_j a_{ij}(x)\partial_j h (x)$$
for all $i, j$ and a.e.-$x$. Then $L$ and $\widetilde L$ are $L^2$-isospectral.}
\dethm

\prf Noting that by the symmetry of the matrix $(a_{ij})$, we have
$$\aligned
L(fh)&=\sum_{i, j}a_{ij}\partial_{ij}^2 (fh)+\sum_i b_i\partial_i (fh)+cfh\\
&=\sum_{i, j}a_{ij}\big[\big(\partial_{ij}^2 f\big) h+ 2\partial_i f \partial_j h\\
&\qd + f\big(\partial_{ij}^2 h\big)\big]
+\sum_i b_i\big[\big(\partial_i f\big) h+f \partial_i h\big] +f ( c h)\\
&=hL f + fL h -cfh +2 \sum_{i, j}a_{ij}\partial_jh \partial_i f
\qd\text{a.e.}
\endaligned
$$
Because $h$ is $L$-harmonic, we obtain
$$\aligned
\frac{1}{h}L(fh)&=(L f- c f) + \frac{2}{h}
\sum_{i}\bigg(\sum_j a_{ij} \partial_j h\bigg) \partial_i f,\qd\text{a.e.}
\endaligned
$$
From which, one reads out the coefficients ${\tilde a}_{ij}(x)$ and ${\tilde b}_i(x)$ of $\widetilde L$.
\deprf

For short, if we set $L_0=L-c$, then we have
$$\aligned
\widetilde L&=L_0+\frac{2}{h}\langle a \nabla h, \nabla \rangle \\
&=L_0+{2}\langle a \nabla \log h, \nabla \rangle\qd \text{if }h>0.
\endaligned$$

\rmk{\cms In one-dimensional case, denoting by $(a(x), b(x), c(x))$ the coefficients of $L$,
we can represent $L$ as
$$L=\frac{\d}{\d\mu}\frac{\d }{\d {\hat\nu}}+ c(x),
$$
where
$$\d\mu(x)=\frac{e^{C(x)}}{a(x)}\d x,\qd \d {\hat\nu}(x)=e^{-C(x)}\d x, \qd C(x)=\int_{\uz}^x \frac{b}{a}(z)\d z,$$
and $\uz$ is a reference point. Then the (dual) operator $\widetilde L$
can be written as
$$\widetilde L=\frac{\d}{\d{\tilde \mu}}\frac{\d }{\d {\hat{\tilde\nu}}}
=\frac{\d}{\d(h^2\mu)}\frac{\d }{\d \big({h^{-2}\hat\nu}\big)}.$$
}
\dermk

Here are simple examples of $L$-harmonic functions.

\xmp{\cms Let $E={\mathbb R}$ or $(0, \infty)$.
\bg{itemize}
\item[(1)] The function $h(x)=x$ is $L$-harmonic (a.e.) on $E$ for
    $$L= \gz(x)\big(\partial_{xx}^2 + V(x)\partial_x - V(x)/x\big),$$
where the functions $V$ and $\gz$ are arbitrary.
\item[(2)] The function $h(x)=x^2$ is $L$-harmonic (a.e.) on $E$ for
    $$L= \gz(x)\big(x\partial_{xx}^2 +\partial_{x} - 4/x\big),$$
where the function $\gz$ is again arbitrary.
\end{itemize}
}
\dexmp

In dimension one, the existence and uniqueness of $L$-harmonic function,
as well as an approximating (constructing) procedure, can be found from
\rf{za05}{Theorems 1.2.1 and 2.2.1}. To see the positivity of $h$ in general dimensions, suppose
that $L$ is self-adjoint and $\sup_x c(x)\le 0$. Then the spectrum of $-L$
should be nonnegative. If the principal eigenvalue $\lz_0$ of $L$ (i.e. the minimal eigenvalue
of $-L$) is zero, then, the $L$-harmonic
function is just a non-trivial eigenfunction corresponding to the eigenvalue $\lz_0=0$
and hence should be nonnegative. The function $h$ should be positive inside
the domain based on the maximum principal. Next, if $\lz_0>0$, then replacing
$L$ by a shift $L+\lz_0$, its principal eigenvalue becomes zero, we can
continue the study as above, and finally shifting back to the original operator.

In higher dimensional case, the harmonic function may not be unique. We
remark that the positive solution of $L$-harmonic functions for Schr\"odinger
operator $L=\ddz +c(x)$ was examined in \ct{mm86} in detail, and for
elliptic operators in \ct{prg95} with probabilistic representation.

\xmp[\rf{mm86}{(1.2)}]{\cms The $L$-harmonic function $h$ for
$L=\ddz -1$ can be represented as
$$h(x)=\int_{S^{n-1}} e^{x\cdot \oz}\d\mu (\oz),$$
where $\mu$ is a nonnegative measure on the unique sphere $S^{n-1}$.
}
\dexmp

The next example is a particular case of Corollary \ref{t1-2}.
Its duality relation was mentioned in \rf{jk12}{\S 6. Example of O.U.-process
and harmonic oscillator}, without mention the $L$-harmonic property of $h$.

\xmp\lb{t3-5}{\cms  On ${\mathbb R}$, the function $h(x)=\exp[-x^2/2]$ is $L$-harmonic:
$$L=\frac 1 2\bigg( \frac{\d^2}{\d x^2 } + 1- x^2\bigg).$$
Its dual is the O.U.-operator:
     $$\widetilde L= \frac 1 2 \frac{\d^2}{\d x^2 } - x\frac{\d}{\d x}.$$
Furthermore, $L$ has $L^2$-eigenvalues $\lz_n=n\,(n\ge 0)$ with eigenfunctions
$$g_n(x)=(-1)^n e^{x^2/2} \frac{\d^n}{\d x^n}\big(e^{-x^2}\big),\qqd n\ge 0,$$
respectively.}
\dexmp

We have just seen an example of the application of known results having
${\tilde c}(x)=0$ to the one having $c(x)\ne 0$. This indicates a
general result as follows.

\thm\lb{t3-6}{\cms Given an elliptic operator
$$\widetilde L= \sum_{i, j}{\tilde a}_{ij}(x)\partial_{ij}^2+ \sum_i {\tilde b}_i(x)\partial_i,
\qqd {\scr D}\big(\tilde L\big)\subset L^2\big(\tilde\mu\big),$$
for each $h\in {\scr C}^2$, $h\ne 0$ a.e., ${\widetilde L}$ is $L^2$-isospectral to $L$:
$$L= \sum_{i, j}{a}_{ij}(x)\partial_{ij}^2+ \sum_i {b}_i(x)\partial_i+ c(x),
\qqd {\scr D}(L)=\big\{f\in {\scr E}: f/h\in {\scr D}\big(\widetilde L\big)\big\},$$
where
$$\aligned
{a}_{ij}(x)&={\tilde a}_{ij}(x),\\
{b}_i(x)&={\tilde b}_i(x)- \frac{2}{h(x)}\sum_j {\tilde a}_{ij}(x) \partial_j h(x) \text{\cms\; on } [h\ne 0],\\
c(x)&=\frac{2}{h(x)^2}\sum_{i, j} {{\tilde a}_{ij}(x)} \partial_i
h(x)\partial_j h(x) -\frac{1}{h(x)}{\widetilde L} h(x) \qd
\text{\cms\; on } [h\ne 0].
\endaligned
$$
Briefly,
$$\aligned
L&=\widetilde L- \frac{2}{h}\big\langle {\tilde a} \nabla h, \nabla \big\rangle
  +\bigg[\frac{2}{h^2}\big\langle {\tilde a} \nabla h, \nabla h \big\rangle
       -\frac{1}{h}{\widetilde L}h \bigg] \\
&=\widetilde L- 2\big\langle {\tilde a} \nabla\log h, \nabla \big\rangle
  + \Big\{2\big\langle {\tilde a} \nabla\log h, \nabla\log h \big\rangle
       -h^{-1}\big\langle {\tilde a} \nabla, \nabla h\big\rangle \\
&\text{\hskip4.2truecm}      + \big\langle {\tilde b}, \nabla\log
h\big\rangle \Big\}\qd \text{\cms if }h>0.
\endaligned$$
}
\dethm

\prf In parallel to the pure jump case, this is simply a use of the duality $L=H {\widetilde L}H^{-1}$,
noting the property that $Lh=0$ is now automatic since $\widetilde L 1=0$.
The remainder of the proof is mainly a careful computation. Actually,
$${\widetilde L}\bigg(\frac f h\bigg)=\frac{1}{h}{\widetilde L} f+ f {\widetilde L}\bigg(\frac{1}{h}\bigg)
+ 2\bigg\langle {\tilde a}\nabla \bigg(\frac{1}{h}\bigg), \nabla f\bigg\rangle.$$
Hence
$$h{\widetilde L}\bigg(\frac f h\bigg)={\widetilde L} f
+2 h\bigg\langle {\tilde a} \nabla \bigg(\frac{1}{h}\bigg), \nabla f\bigg\rangle
+ fh {\widetilde L}\bigg(\frac{1}{h}\bigg).$$
From this, it is ready to write down the coefficients of $L$.\deprf

\crl\lb{t3-7}{\cms For given $\widetilde L$ and $h=\exp \psi$, the dual operator $L$ takes the following
form
$$L=\widetilde L - 2\big\langle {\tilde a} \nabla \psi, \nabla \big\rangle
  + \Big\{\big\langle {\tilde a} \nabla \psi, \nabla \psi \big\rangle
            -{\widetilde L} \psi\Big\}.$$
}
\decrl

We remark that Corollary \ref{t3-7} provides us an alternative way to construct the
isospectral operator in dimension one. Suppose that we are given an operator
$${\overline L}= {\bar a}(x)\frac{\d^2}{\d x^2} + {\bar b}(x)\frac{\d}{\d x} + {\bar c}(x).$$
We want to construct $\widetilde L$ in terms of the operator $L$ given in Corollary \ref{t3-7}.
First, instead of solving the second order harmonic equation $\overline L h=0$, we need
to solve the first order Riccati equation for $\phi$:
$${\bar a}\phi'+{\bar a}\phi^2+{\bar b}\phi +{\bar c}=0$$
to which there is a standard iterative procedure in ODE. Next, let $\psi$ satisfy $\psi'=\phi$
and define $\tilde b= 2{\bar a}\phi+{\bar b} $. Then we have $L={\overline L}$. With this
$\tilde b$ and $\tilde a:={\bar a}$, we obtain the operator $\widetilde L$ as required.

As an application of the last theorem, one can obtain a lot of
examples from \ct{cmf12a, cmf12b}. We remark that each $\widetilde
L$ corresponds to a large class of $L$ since $h$ is quite
arbitrary.

The natural higher-dimensional extension of Example \ref{t3-5} is as follows.

\xmp{\cms The dual of
    $L=\frac 1 2\sum_i \big(\partial_{ii}^2 + 1- x_i^2\big)$
is
     $\widetilde L= \frac 1 2 \sum_i\big(\partial_{ii}^2 -2  x_i\partial_i \big).$
The function $h$ takes the form
   $h(x)=\exp[-|x|^2/2]$
rather than
   $\sum_i \exp\big[-x_i^2/2\big].$
The operator $L$ has eigenvalue $n\,(n\ge 0)$ with multiplicity
$\#\{(k_1, k_2, \ldots, k_d): k_1+ k_2+ \ldots+ k_d=n\}$, here
$\#$ means the cardinality of the set following. } \dexmp

\prf For the higher-dimensional O.U.-operator $\widetilde L$, we have
eigenvalues $\{\sum_{i=1}^d k_i: k_i=0,1,\ldots\}$. Corresponding to each
$\sum_{i=1}^d k_i$, the eigenfunction is $g(x):=\prod_{i=1}^d g_{k_i}^{(i)}(x_i)$\,
(where each $g_n^{(i)}$ is the function $g_n$ given in the proof of Corollary \ref{t1-2}):
$$\widetilde L g(x)= - \sum_{i=1}^d k_i g_{k_i}^{(i)}(x_i)\prod_{j\ne i} g_{k_j}^{(j)}(x_j)
=-\bigg(\sum_{i=1}^d k_i\bigg) g(x).$$
Therefore, $\widetilde L$ has eigenvalue $n\,(n\ge 0)$ with multiplicity
$\#\{(k_1, k_2, \ldots, k_d): k_1+ k_2+ \ldots+ k_d=n\}$.
From here, it is easy to write down the eigenvalues of $L$ and their corresponding
eigenfunctions.
\deprf

\medskip

\nnd{\bf Acknowledgments}. {\small The results of
the paper were presented several times in our seminar, from which the authors are benefited a lot
from the discussions and suggestions. Research supported in part by the
         National Natural Science Foundation of China (No. 11131003),
         the ``985'' project from the Ministry of Education in China,
and the Project Funded by the Priority Academic Program Development of Jiangsu Higher Education Institutions.
}

\medskip

\nnd {\small
Mu-Fa Chen\\
School of Mathematical Sciences, Beijing Normal University,
Laboratory of Mathematics and Complex Systems (Beijing Normal University),
Ministry of Education, Beijing 100875,
    The People's Republic of China.\newline E-mail: mfchen@bnu.edu.cn\newline Home page:
    http://math.bnu.edu.cn/\~{}chenmf/main$\_$eng.htm\\
{}\\
Xu Zhang\\
College of Applied Sciences, Beijing University of Technology,
Beijing 100022,
    The People's Republic of China.
    \newline E-mail: zhangxu@bjut.edu.cn}

\end{document}